# Distribution of the partition function modulo $m$

By KEN ONO*

## 1. Introduction and statement of results

A *partition* of a positive integer $n$ is any nonincreasing sequence of positive integers whose sum is $n$. Let $p(n)$ denote the number of partitions of $n$ (as usual, we adopt the convention that $p(0) = 1$ and $p(\alpha) = 0$ if $\alpha \notin \mathbb{N}$). Ramanujan proved for every nonnegative integer $n$ that

$$p(5n + 4) \equiv 0 \pmod{5},$$
$$p(7n + 5) \equiv 0 \pmod{7},$$
$$p(11n + 6) \equiv 0 \pmod{11},$$

and he conjectured further such congruences modulo arbitrary powers of 5, 7, and 11. Although the work of A. O. L. Atkin and G. N. Watson settled these conjectures many years ago, the congruences have continued to attract much attention. For example, subsequent works by G. Andrews, A. O. L. Atkin, F. Garvan, D. Kim, D. Stanton, and H. P. F. Swinnerton-Dyer ([An-G], [G], [G-K-S], [At-Sw2]), in the spirit of F. Dyson, have gone a long way towards providing combinatorial and physical explanations for their existence.

Ramanujan [Ra, p. xix] already observed that his congruences were quite special. For instance, he proclaimed that

"*It appears that there are no equally simple properties for any moduli involving primes other than these three (i.e. $m = 5, 7, 11$).*"

Although there is no question that congruences of the form $p(an + b) \equiv 0 \pmod{m}$ are rare (see recent works by the author ([K-Ol], [O1], [O2])), the question of whether there are many such congruences has been the subject of debate. In the 1960's, Atkin and O'Brien ([At], [At-Sw1], [At-Ob]) uncovered

---

*The author is supported by NSF grants DMS-9508976, DMS-9874947 and NSA grant MSPR-97Y012.

1991 *Mathematics Subject Classification*. Primary 11P83; Secondary 05A17.

*Key words and phrases.* partition function, The Erdös' conjecture, Newman's conjecture.

further congruences such as

(1) $$p(11^3 \cdot 13n + 237) \equiv 0 \pmod{13}.$$

However, no further congruences have been found and proven since.

In a related direction, P. Erdös and A. Ivić ([E-I]) conjectured that there are infinitely many primes $m$ which divide some value of the partition function, and Erdös made the following stronger conjecture [Go], [I].

CONJECTURE (Erdös).   *If $m$ is prime, then there is at least one nonnegative integer $n_m$ for which*

$$p(n_m) \equiv 0 \pmod{m}.$$

A. Schinzel (see [E-I] for the proof) proved the Erdös-Ivić conjecture using the Hardy-Ramanujan-Rademacher asymptotic formula for $p(n)$, and more recently Schinzel and E. Wirsing [Sc-W] have obtained a quantitative result in the direction of Erdös' stronger conjecture. They have shown that the number of primes $m < X$ for which Erdös' conjecture is true is $\gg \log\log X$.

Here we present a uniform and systematic approach which settles the debate regarding the existence of further congruences, and yields Erdös' conjecture as an immediate corollary.

THEOREM 1.   *Let $m \geq 5$ be prime and let $k$ be a positive integer. A positive proportion of the primes $\ell$ have the property that*

$$p\left(\frac{m^k \ell^3 n + 1}{24}\right) \equiv 0 \pmod{m}$$

*for every nonnegative integer $n$ coprime to $\ell$.*

In view of work of S. Ahlgren [A], J.-L. Nicolas, I. Z. Ruzsa, A. Sárközy [Ni-R-Sa] and J-P. Serre [S] for $m = 2$, the fact that $p(3) = 3$, and Theorem 1, we obtain:

COROLLARY 2.   *Erdös' conjecture is true for every prime $m$. Moreover, if $m \neq 3$ is prime, then*

$$\#\{0 \leq n \leq X \ : \ p(n) \equiv 0 \pmod{m}\} \gg_m \begin{cases} \sqrt{X} & \text{if } m = 2, \\ X & \text{if } m \geq 5. \end{cases}$$

Surprisingly, it is not known whether there are infinitely many $n$ for which $p(n) \equiv 0 \pmod{3}$.

As an example, we shall see that $\ell = 59$ satisfies the conclusion of Theorem 1 when $m = 13$ and $k = 1$. In this case, by considering integers in the arithmetic progression $r \equiv 1 \pmod{24 \cdot 59}$, we find for every nonnegative integer $n$ that

(2) $$p(59^4 \cdot 13n + 111247) \equiv 0 \pmod{13}.$$

Our results are also useful in attacking a famous conjecture of M. Newman [N1].

CONJECTURE (M. Newman). *If $m$ is an integer, then for every residue class $r \pmod{m}$ there are infinitely many nonnegative integers $n$ for which $p(n) \equiv r \pmod{m}$.*

Works by Atkin, Newman, and O. Kolberg ([At], [N1], [K]) have verified the conjecture for $m = 2, 5, 7, 11$ and $13$ (in fact, the case where $m = 11$ is not proved in these papers, but one may easily modify the arguments to obtain this case). Here we present a result which, in principle, may be used to verify Newman's conjecture for every remaining prime $m \neq 3$.

We shall call a prime $m \geq 5$ *good* if for every $r \pmod{m}$ there is a nonnegative integer $n_r$ for which $mn_r \equiv -1 \pmod{24}$ and

$$p\left(\frac{mn_r + 1}{24}\right) \equiv r \pmod{m}.$$

THEOREM 3. *If $m \geq 5$ is a good prime, then Newman's conjecture is true for $m$. Moreover, for each residue class $r \pmod{m}$ we have*

$$\#\{0 \leq n \leq X \ : \ p(n) \equiv r \pmod{m}\} \gg_{r,m} \begin{cases} \sqrt{X}/\log X & \text{if } 1 \leq r \leq m - 1, \\ X & \text{if } r = 0. \end{cases}$$

Although it appears likely that every prime $m \geq 13$ is good, proving that a prime $m$ is good involves a substantial computation, and this computation becomes rapidly infeasible as the size of $m$ grows. The author is indebted to J. Haglund and C. Haynal who wrote efficient computer code to attack this problem. As a result, we have the following.

COROLLARY 4. *Newman's conjecture is true for every prime $m < 1000$ with the possible exception of $m = 3$.*

We also uncover surprising "periodic" relations for certain values of the partition function mod $m$. In particular, we prove that if $m \geq 5$ is prime, then the sequence of generating functions

(3) $$F(m, k; z) := \sum_{\substack{n \geq 0 \\ m^k n \equiv -1 \pmod{24}}} p\left(\frac{m^k n + 1}{24}\right) q^n \pmod{m}$$

($q := e^{2\pi i z}$ throughout) is eventually periodic in $k$. We call these periods "Ramanujan cycles." Their existence implies the next result.

THEOREM 5. *If $m \geq 5$ is prime, then there are integers $0 \leq N(m) \leq 48(m^3 - 2m - 1)$ and $1 \leq P(m) \leq 48(m^3 - 2m - 1)$ such that for every $i > N(m)$ we have*

$$p\left(\frac{m^i n + 1}{24}\right) \equiv p\left(\frac{m^{P(m)+i} \cdot n + 1}{24}\right) \pmod{m}$$

*for every nonnegative integer $n$.*

For each class $r \pmod{m}$ one obtains explicit sequences of integers $n_k$ such that $p(n_k) \equiv r \pmod{m}$ for all $k$. This is the subject of Corollaries 9 through 12 below. For example, taking $n = 0$ in Corollary 12 shows that for every nonnegative integer $k$

(4)
$$p\left(\frac{23^{2k+1} + 1}{24}\right) \equiv 5^k \pmod{23} \quad \text{and} \quad p\left(\frac{23^{2k+3} + 1}{24}\right) \equiv 5^{k+1} \pmod{23}.$$

Similarly it is easy to show that

(5) $$p\left(\frac{1367 \cdot 23^{2k+2} + 1}{24}\right) \equiv 0 \pmod{23} \quad \text{and}$$

$$p\left(\frac{1297 \cdot 23^{2k+1} + 1}{24}\right) \equiv 0 \pmod{23}.$$

Congruences of this sort mod 13 were previously discovered by Ramanujan and found by M. Newman [N2]. In fact, this paper was inspired by such entries in Ramanujan's lost manuscript on $p(n)$ and $\tau(n)$ (see [B-O]).

*A priori*, one knows that the generating functions $F(m, k; z)$ are the reductions mod $m$ of weight $-1/2$ nonholomorphic modular forms, and as such lie in infinite dimensional $\mathbb{F}_m$-vector spaces. This infinitude has been the main obstacle in obtaining results for the partition function mod $m$. In Section 3 we shall prove a theorem (see Theorem 8) which establishes that the $F(m, k; z)$ are the reductions mod $m$ of half-integral weight cusp forms lying in one of two spaces with Nebentypus. Hence, there are only finitely many possibilities for each $F(m, k; z)$. This is the main observation which underlies all of the results in this paper. We then prove Theorems 1 and 3 by employing the Shimura correspondence and a theorem of Serre about Galois representations. In Section 4 we present detailed examples for $5 \leq m \leq 23$.

## 2. Preliminaries

We begin by defining operators $U$ and $V$ which act on formal power series. If $M$ and $j$ are positive integers, then

$$(6) \qquad \left(\sum_{n\geq 0} a(n)q^n\right) \mid U(M) := \sum_{n\geq 0} a(Mn)q^n,$$

$$(7) \qquad \left(\sum_{n\geq 0} a(n)q^n\right) \mid V(j) := \sum_{n\geq 0} a(n)q^{jn}.$$

We recall that Dedekind's eta-function is defined by

$$(8) \qquad \eta(z) := q^{1/24} \prod_{n=1}^{\infty}(1-q^n)$$

and that Ramanujan's Delta-function is

$$(9) \qquad \Delta(z) := \eta^{24}(z),$$

the unique normalized weight 12 cusp form for $\mathrm{SL}_2(\mathbb{Z})$. If $m \geq 5$ is prime and $k$ is a positive integer, then define $a(m,k,n)$ by

$$(10) \qquad \sum_{n=0}^{\infty} a(m,k,n)q^n := \frac{\left(\Delta^{\delta(m,k)}(z) \mid U(m^k)\right) \mid V(24)}{\eta^{m^k}(24z)} \pmod{m},$$

where $\delta(m,k) := (m^{2k}-1)/24$. Recall the definition (3) of $F(m,k;z)$.

THEOREM 6. *If $m \geq 5$ is prime and $k$ is a positive integer, then*

$$F(m,k;z) \equiv \sum_{n=0}^{\infty} a(m,k,n)q^n \pmod{m}.$$

*Proof.* We begin by recalling that Euler's generating function for $p(n)$ is given by the infinite product

$$\sum_{n=0}^{\infty} p(n)q^n := \prod_{n=1}^{\infty} \frac{1}{(1-q^n)}.$$

Using this fact, one easily finds that

$$\frac{\eta^{m^k}(m^k z)}{\eta(z)} \mid U(m^k) = \left\{\sum_{n=0}^{\infty} p(n)q^{n+\delta(m,k)} \cdot \prod_{n=1}^{\infty}(1-q^{m^k n})^{m^k}\right\} \mid U(m^k)$$

$$= \sum_{n=0}^{\infty} p(m^k n + \beta(m,k))q^{n+\frac{\delta(m,k)+\beta(m,k)}{m^k}} \cdot \prod_{n=1}^{\infty}(1-q^n)^{m^k},$$

where $1 \leq \beta(m,k) \leq m^k - 1$ satisfies $24\beta(m,k) \equiv 1 \pmod{m^k}$.

Since $(1 - X^{m^k})^{m^k} \equiv (1 - X)^{m^{2k}} \pmod{m}$, we find that

$$\sum_{n=0}^{\infty} p(m^k n + \beta(m,k)) q^{n + \frac{\delta(m,k) + \beta(m,k)}{m^k}} \equiv \frac{\Delta^{\delta(m,k)}(z) \mid U(m^k)}{\prod_{n=1}^{\infty}(1-q^n)^{m^k}} \pmod{m}.$$

Replacing $q$ by $q^{24}$ and multiplying through by $q^{-m^k}$ one obtains

$$\sum_{n=0}^{\infty} p(m^k n + \beta(m,k)) q^{24n + \frac{24\beta(m,k)-1}{m^k}} \equiv \sum_{n=0}^{\infty} a(m,k,n) q^n \pmod{m}.$$

It is easy to see that

$$\sum_{n=0}^{\infty} p(m^k n + \beta(m,k)) q^{24n + \frac{24\beta(m,k)-1}{m^k}} = \sum_{\substack{n \geq 0 \\ m^k n \equiv -1 \pmod{24}}} p\left(\frac{m^k n + 1}{24}\right) q^n. \quad \square$$

We conclude this section with the following elementary result which establishes that the $F(m,k;z)$ form an inductive sequence generated by the action of the $U(m)$ operator.

PROPOSITION 7. *If $m \geq 5$ is prime and $k$ is a positive integer, then*

$$F(m, k+1; z) \equiv F(m, k; z) \mid U(m) \pmod{m}.$$

*Proof.* Using definition of the $F(m,k;z)$ and the convention that $p(\alpha) = 0$ for $\alpha \notin \mathbb{Z}$, one finds that

$$F(m,k;z) \mid U(m) \equiv \sum_{\substack{n \geq 0 \\ m^k n \equiv -1 \pmod{24}}} p\left(\frac{m^k n + 1}{24}\right) q^n \mid U(m)$$

$$= \sum_{\substack{n \geq 0 \\ m^{k+1} n \equiv -1 \pmod{24}}} p\left(\frac{m^{k+1} n + 1}{24}\right) q^n$$

$$\equiv F(m, k+1; z) \pmod{m}. \quad \square$$

## 3. Proof of the results

First we recall some notation. Suppose that $w \in \frac{1}{2}\mathbb{Z}$, and that $N$ is a positive integer (with $4 \mid N$ if $w \notin \mathbb{Z}$). Let $S_w(\Gamma_0(N), \chi)$ denote the space of weight $w$ cusp forms with respect to the congruence subgroup $\Gamma_0(N)$ and with Nebentypus character $\chi$. Moreover, if $\ell$ is prime, then let $S_w(\Gamma_0(N), \chi)_\ell$ denote the $\mathbb{F}_\ell$-vector space of the reductions mod $\ell$ of the $q$-expansions of forms in $S_w(\Gamma_0(N), \chi)$ with rational integer coefficients.

THEOREM 8. *If $m \geq 5$ is prime, then for every positive integer $k$ we have*

$$F(m, k; z) \in S_{\frac{m^2-m-1}{2}}(\Gamma_0(576m), \chi\chi_m^{k-1})_m,$$

*where $\chi$ is the nontrivial quadratic character with conductor 12, and $\chi_m$ is the usual Kronecker character for $\mathbb{Q}(\sqrt{m})$.*

*Proof.* The $U(m)$ operator defines a map (see [S-St, Lemma 1])

$$U(m): S_{\lambda+\frac{1}{2}}(\Gamma_0(4Nm), \nu) \longrightarrow S_{\lambda+\frac{1}{2}}(\Gamma_0(4Nm), \nu\chi_m).$$

Therefore, in view of Proposition 7 it suffices to prove that

$$F(m, 1; z) \in S_{\frac{m^2-m-1}{2}}(\Gamma_0(576m), \chi)_m.$$

If $d \equiv 0 \pmod{4}$, then it is well known that the space of cusp forms $S_d(\Gamma_0(1))$ has a basis of the form

$$\left\{ \Delta(z)^j E_4(z)^{\frac{d}{4}-3j} \; : \; 1 \leq j \leq \left[\frac{d}{12}\right] \right\}.$$

Since the Hecke operator $T_m$ is the same as the $U(m)$ operator on $S_{12\delta(m,1)}(\Gamma_0(1))_m$, we know that

$$\Delta^{\delta(m,1)}(z) \mid U(m) \equiv \sum_{j \geq 1} \alpha_j \Delta(z)^j E_4(z)^{3\delta(m,1)-3j} \pmod{m},$$

where the $\alpha_j \in \mathbb{F}_m$. However, since

$$\Delta^{\delta(m,1)}(z) = q^{\delta(m,1)} - \cdots,$$

it is easy to see that

$$\Delta^{\delta(m,1)}(z) \mid U(m) = \sum_{n \geq n_0} t(n) q^n$$

where $n_0 \geq \delta(m,1)/m$. However, since $\delta(m,1) \in \mathbb{Z}$, one can easily deduce that $n_0 > m/24$.

The only basis forms in $\Delta^{\delta(m,1)}(z) \mid U(m) \pmod{m}$ are those $\Delta^j(z) E_4(z)^{3\delta(m,1)-3j}$ where $j > m/24$. This implies that

$$\frac{\left(\Delta(z)^{\delta(m,1)} \mid U(m)\right) \mid V(24)}{\eta^m(24z)}$$

is a cusp form. Since $\left(\Delta(z)^{\delta(m,1)} \mid U(m)\right) \mid V(24)$ is the reduction mod $m$ of a weight $\frac{m^2-1}{2}$ cusp form with respect to $\Gamma_0(24)$, and $\eta(24z)$ is a weight $1/2$ cusp form with respect to $\Gamma_0(576)$ with character $\chi$, the result follows. □

Now we recall an important result due to Serre [S, 6.4].

THEOREM (Serre). *The set of primes $\ell \equiv -1 \pmod{N}$ for which*

$$f \mid T_\ell \equiv 0 \pmod{m}$$

*for every $f(z) \in S_k(\Gamma_0(N), \nu)_m$ has positive density. Here $T_\ell$ denotes the usual Hecke operator of index $\ell$ acting on $S_k(\Gamma_0(N), \nu)$.*

*Proof of Theorem* 1. If $F(m, k; z) \equiv 0 \pmod{m}$, then the conclusion of Theorem 1 holds for every prime $\ell$. Hence, we may assume that $F(m, k; z) \not\equiv 0 \pmod{m}$. By Theorem 8, we know that each $F(m, k; z)$ belongs to $S_{\frac{m^2-m-1}{2}}(\Gamma_0(576m), \chi \chi_m^{k-1})_m$. Therefore each $F(m, k; z)$ is the reduction mod $m$ of a half-integral weight cusp form.

Now we briefly recall essential facts about the "Shimura correspondence" ([Sh]), a family of maps which send modular of forms of half-integral weight to those of integer weight. Although Shimura's original theorem was stated for half-integral weight eigenforms, the generalization we describe here follows from subsequent works by Cipra and Niwa [Ci], [Ni]. Suppose that $f(z) = \sum_{n=1}^\infty b(n) q^n \in S_{\lambda+\frac{1}{2}}(\Gamma_0(4N), \psi)$ is a cusp form where $\lambda \geq 2$. If $t$ is any square-free integer, then define $A_t(n)$ by

$$\sum_{n=1}^\infty \frac{A_t(n)}{n^s} := L(s - \lambda + 1, \psi \chi_{-1}^\lambda \chi_t) \cdot \sum_{n=1}^\infty \frac{b(tn^2)}{n^s}.$$

Here $\chi_{-1}$ (resp. $\chi_t$) is the Kronecker character for $\mathbb{Q}(i)$ (resp. $\mathbb{Q}(\sqrt{t})$). These numbers $A_t(n)$ define the Fourier expansion of $S_t(f(z))$, a cusp form

$$S_t(f(z)) := \sum_{n=1}^\infty A_t(n) q^n$$

in $S_{2\lambda}(\Gamma_0(4N), \psi^2)$. Moreover, the Shimura correspondence $S_t$ commutes with the Hecke algebra. In other words, if $p \nmid 4N$ is prime, then

$$S_t(f \mid T(p^2)) = S_t(f) \mid T_p.$$

Here $T_p$ (resp. $T(p^2)$) denotes the usual Hecke operator acting on the space $S_{2\lambda}(\Gamma_0(4N), \psi^2)$ (resp. $S_{\lambda+\frac{1}{2}}(\Gamma_0(4N), \psi)$).

Therefore, for every square-free integer $t$ we have that the image $S_t(F(m, k; z))$ under the $t^{\text{th}}$ Shimura correspondence is the reduction mod $m$ of an integer weight form in $S_{m^2-m-2}(\Gamma_0(576m), \chi_{\text{triv}})$. Now let $S(m)$ denote the set of primes $\ell \equiv -1 \pmod{576m}$ for which

$$G \mid T_\ell \equiv 0 \pmod{m}$$

for every $G \in S_{m^2-m-2}(\Gamma_0(576m), \chi_{\text{triv}})_m$. By Serre's theorem, the set $S(m)$ contains a positive proportion of the primes.

By the commutativity of the correspondence, if $\ell \in S(m)$, then we find that
$$F(m, k; z) \mid T(\ell^2) \equiv 0 \pmod{m},$$
where $T(\ell^2)$ is the Hecke operator of index $\ell^2$ on $S_{\frac{m^2-m-1}{2}}(\Gamma_0(576m), \chi\chi_m^{k-1})$.
In particular (see [Sh]), if $f = \sum a_f(n)q^n \in S_{\lambda+\frac{1}{2}}(N, \chi_f)$ is a half-integral weight form, then

$$(11) \quad f \mid T(\ell^2) := \sum_{n=0}^{\infty} \left( a_f(\ell^2 n) + \chi_f(\ell)\left(\frac{(-1)^\lambda n}{\ell}\right)\ell^{\lambda-1} a_f(n) \right.$$
$$\left. + \chi_f(\ell^2)\ell^{2\lambda-1} a_f(n/\ell^2) \right) q^n.$$

Therefore, if $\ell \in S(m)$ and $n$ is a positive integer which is coprime to $\ell$, then, by replacing $n$ by $n\ell$, we have

$$a(m, k, n\ell^3) + \chi\chi_m^{k-1}(\ell)\left(\frac{(-1)^{\frac{m^2-m-2}{2}} n\ell}{\ell}\right) \cdot \ell^{\frac{m^2-m-4}{2}} \cdot a(m, k, n\ell) \equiv 0 \pmod{m}.$$

Since $\left(\frac{n\ell}{\ell}\right) = 0$, by Theorem 6 we find that
$$p\left(\frac{m^k \ell^3 n + 1}{24}\right) \equiv a(m, k, \ell^3 n) \equiv 0 \pmod{m}. \qquad \square$$

*Remark.* Although Theorem 1 is a general result guaranteeing the existence of congruences, there are other congruences which follow from other similar arguments based on (11).

For example, suppose that $\ell$ is a prime for which
$$F(m, k; z) \mid T(\ell^2) \equiv \lambda(\ell) F(m, k; z) \pmod{m}$$
for some $\lambda(\ell) \in \mathbb{F}_m^\times$. If $n$ is a nonnegative integer for which $\ell^2 \nmid n$, then (11) becomes

$$a(m, k, n)\left\{ \lambda(\ell) - \chi\chi_m^{k-1}(\ell)\left(\frac{(-1)^{\frac{m^2-m-2}{2}} n}{\ell}\right)\ell^{\frac{m^2-m-4}{2}} \right\}$$
$$\equiv a(m, k, n\ell^2) \pmod{m}.$$

Hence, if it turns out that
$$\lambda(\ell) \equiv \pm \ell^{\frac{m^2-m-4}{2}} \pmod{m},$$
then there are arithmetic progressions of integers $n$ for which
$$a(m, k, n\ell^2) \equiv p\left(\frac{m^k \ell^2 n + 1}{24}\right) \equiv 0 \pmod{m}.$$

Although we have not conducted a thorough search, it is almost certain that many such congruences exist.

*Proof of Theorem* 3. By the proof of Theorem 6, recall that

$$F(m,1;z) = \sum_{\substack{n \geq 0, \\ mn \equiv -1 \pmod{24}}} p\left(\frac{mn+1}{24}\right) q^n \in S_{\frac{m^2-m-1}{2}}(\Gamma_0(576m), \chi)_m.$$

Since $m$ is good, for each $0 \leq r \leq m-1$ let $n_r$ be a fixed nonnegative integer for which $mn_r \equiv -1 \pmod{24}$ and

$$p\left(\frac{mn_r+1}{24}\right) \equiv r \pmod{m}.$$

Let $M_m$ be the set of primes $p$ for which $p \mid n_r$ for some $r$, and define $\mathfrak{S}_m$ by

$$\mathfrak{S}_m := \prod_{p \in M_m} p.$$

Obviously, the form $F(m,1;z)$ also lies in $S_{\frac{m^2-m-1}{2}}(\Gamma_0(576m\mathfrak{S}_m), \chi)_m$. Therefore, by Serre's theorem and the commutativity of the Shimura correspondence, a positive proportion of the primes $\ell \equiv -1 \pmod{576m\mathfrak{S}_m}$ have the property that

$$F(m,1;z) \mid T(\ell^2) \equiv 0 \pmod{m}.$$

By (11), for all but finitely many such $\ell$ we have for each $r$ that

$$p\left(\frac{mn_r\ell^2+1}{24}\right) + \chi(\ell)\left(\frac{(-1)^{\frac{m^2-m-2}{2}} n_r}{\ell}\right) \ell^{\frac{m^2-m-4}{2}} p\left(\frac{mn_r+1}{24}\right) \equiv 0 \pmod{m}.$$

However, since $\ell \equiv -1 \pmod{m}$ this implies that

$$(12) \quad p\left(\frac{mn_r\ell^2+1}{24}\right) \equiv \chi(\ell)\left(\frac{(-1)^{\frac{m^2-m-2}{2}}}{\ell}\right)(-1)^{\frac{m^2-m-2}{2}} \left(\frac{n_r}{\ell}\right) r \pmod{m}.$$

If $n_r = \prod_i p_i$ where the $p_i$ are prime, then

$$\left(\frac{n_r}{\ell}\right) := \prod_i \left(\frac{p_i}{\ell}\right).$$

Since $n_r$ is odd, $\ell \equiv 3 \pmod{4}$, and $\ell \equiv -1 \pmod{p_i}$, we find by quadratic reciprocity that

$$\left(\frac{p_i}{\ell}\right) = \left(\frac{\ell}{p_i}\right)\left(\frac{-1}{p_i}\right)$$
$$= \left(\frac{-\ell}{p_i}\right) = \left(\frac{1}{p_i}\right) = 1.$$

Therefore, for all but finitely many such $\ell$ congruence (12) reduces to

$$(13) \qquad p\left(\frac{mn_r\ell^2 + 1}{24}\right) \equiv \chi(\ell)\left(\frac{(-1)^{\frac{m^2-m-2}{2}}}{\ell}\right)(-1)^{\frac{m^2-m-2}{2}}r \pmod{m}.$$

Hence for every sufficiently large such $\ell$, the $m$ values $p\left(\frac{mn_r\ell^2+1}{24}\right)$ are distinct and represent each residue class mod $m$.

To complete the proof, it suffices to notice that the number of such primes $\ell < X$, by Serre's theorem again, is $\gg X/\log X$. In view of (13), this immediately yields the $\sqrt{X}/\log X$ estimate. The estimate when $r = 0$ follows easily from Theorem 1. □

*Proof of Theorem 5.* Since $F(m, k; z)$ is in $S_{\frac{m^2-m-1}{2}}(\Gamma_0(576m), \chi\chi_m^{k-1})_m$, it follows that each $F(m, k; z)$ lies in one of two finite-dimensional $\mathbb{F}_m$-vector spaces. The result now follows immediately from (3), Theorem 6, Proposition 7, and well-known upper bounds for the dimensions of spaces of cusp forms (see [C-O]). □

## 4. Examples

In this section we list the Ramanujan cycles for the generating functions $F(m, k; z)$ when $5 \leq m \leq 23$. Although we have proven that each $F(m, k; z) \in S_{\frac{m^2-m-1}{2}}(\Gamma_0(576m), \chi\chi_m^{k-1})_m$, in these examples it turns out that they all are congruent mod $m$ to forms of smaller weight.

*Cases where $m = 5, 7,$ and $11$.* In view of the Ramanujan congruences mod $5, 7,$ and $11$, it is immediate that for every positive integer $k$ we have

$$F(5, k; z) \equiv 0 \pmod 5,$$
$$F(7, k; z) \equiv 0 \pmod 7,$$
$$F(11, k; z) \equiv 0 \pmod{11}.$$

Therefore, these Ramanujan cycles are degenerate.

*Case where $m = 13$.* By [Gr-O, Prop. 4] it is known that

$$\Delta^7(z) \mid U(13) \equiv 11\Delta(z) \pmod{13}.$$

Therefore by (10) and Theorem 6 it turns out that

$$F(13, 1; z) \equiv 11q^{11} + 9q^{35} + \cdots \equiv 11\eta^{11}(24z) \pmod{13}.$$

Using a theorem of Sturm [St, Th. 1], one easily verifies with a finite computation that

$$F(13, 1; z) \mid T(59^2) \equiv 0 \pmod{13}.$$

By the proof of Theorem 1, we find that every nonnegative integer $n \equiv 1 \pmod{24}$ that is coprime to 59 has the property that

$$p\left(\frac{13 \cdot 59^3 n + 1}{24}\right) \equiv 0 \pmod{13}.$$

Congruence (2) follows immediately.

Using Sturm's theorem again, one readily verifies that

$$\eta^{11}(24z) \mid U(13) \equiv 8\eta^{23}(24z) \pmod{13},$$
$$\eta^{23}(24z) \mid U(13) \equiv 4\eta^{11}(24z) \pmod{13}.$$

By Proposition 7 this implies that

$$F(13, 2; z) \equiv 10\eta^{23}(24z) \pmod{13},$$

and more generally it implies that for every nonnegative integer $k$

(14) $\qquad F(13, 2k+1; z) \equiv 11 \cdot 6^k \eta^{11}(24z) \pmod{13},$

(15) $\qquad F(13, 2k+2; z) \equiv 10 \cdot 6^k \eta^{23}(24z) \pmod{13}.$

These two congruences appear in Ramanujan's unpublished manuscript on $\tau(n)$ and $p(n)$, and their presence in large part inspired this entire work. From (14) and (15) we obtain the following easy corollary.

COROLLARY 9. *Define integers $a(n)$ and $b(n)$ by*

$$\sum_{n=0}^{\infty} a(n)q^n := \prod_{n=1}^{\infty}(1-q^n)^{11},$$
$$\sum_{n=0}^{\infty} b(n)q^n := \prod_{n=1}^{\infty}(1-q^n)^{23}.$$

*If $k$ and $n$ are nonnegative integers, then*

$$p\left(\frac{13^{2k+1}(24n+11)+1}{24}\right) \equiv 11 \cdot 6^k \cdot a(n) \pmod{13},$$
$$p\left(\frac{13^{2k+2}(24n+23)+1}{24}\right) \equiv 10 \cdot 6^k \cdot b(n) \pmod{13}.$$

*Case where $m = 17$.* By [Gr-O, Prop. 4], it is known that

$$\Delta^{12}(z) \mid U(17) \equiv 7E_4(z)\Delta(z) \pmod{17}$$

where $E_4(z) = 1 + 240\sum_{n=1}^{\infty} \sigma_3(n)q^n$ is the usual weight 4 Eisenstein series. Therefore by (10) it turns out that

$$F(17, 1; z) \equiv 7q^7 + 16q^{31} + \cdots \equiv 7\eta^7(24z)E_4(24z) \pmod{17}.$$

Again using Sturm's theorem one easily verifies that
$$\eta^7(24z)E_4(24z) \mid U(17) \equiv 7\eta^{23}(24z)E_4(24z) \pmod{17},$$
$$\eta^{23}(24z)E_4(24z) \mid U(17) \equiv 13\eta^7(24z)E_4(24z) \pmod{17}.$$
By Proposition 7 this implies that for every nonnegative integer $k$
$$F(17, 2k+1; z) \equiv 7 \cdot 6^k \eta^7(24z)E_4(24z) \pmod{17},$$
$$F(17, 2k+2; z) \equiv 15 \cdot 6^k \eta^{23}(24z)E_4(24z) \pmod{17}.$$
As an immediate corollary we obtain:

COROLLARY 10. *Define integers $c(n)$ and $d(n)$ by*
$$\sum_{n=0}^{\infty} c(n)q^n := E_4(z) \cdot \prod_{n=1}^{\infty}(1-q^n)^7,$$
$$\sum_{n=0}^{\infty} d(n)q^n := E_4(z) \cdot \prod_{n=1}^{\infty}(1-q^n)^{23}.$$
*If $k$ and $n$ are nonnegative integers, then*
$$p\left(\frac{17^{2k+1}(24n+7)+1}{24}\right) \equiv 7 \cdot 6^k \cdot c(n) \pmod{17},$$
$$p\left(\frac{17^{2k+2}(24n+23)+1}{24}\right) \equiv 15 \cdot 6^k \cdot d(n) \pmod{17}.$$

*Case where $m = 19$.* Using [Gr-O, Prop. 4], and arguing as above it turns out that for every nonnegative integer $k$
$$F(19, 2k+1; z) \equiv 5 \cdot 10^k \eta^5(24z)E_6(24z) \pmod{19},$$
$$F(19, 2k+2; z) \equiv 11 \cdot 10^k \eta^{23}(24z)E_6(24z) \pmod{19}.$$
Here $E_6(z) = 1 - 504\sum_{n=1}^{\infty}\sigma_5(n)q^n$ is the usual weight 6 Eisenstein series. As an immediate corollary we obtain:

COROLLARY 11. *Define integers $e(n)$ and $f(n)$ by*
$$\sum_{n=0}^{\infty} e(n)q^n := E_6(z) \cdot \prod_{n=1}^{\infty}(1-q^n)^5,$$
$$\sum_{n=0}^{\infty} f(n)q^n := E_6(z) \cdot \prod_{n=1}^{\infty}(1-q^n)^{23}.$$
*If $k$ and $n$ are nonnegative integers, then*
$$p\left(\frac{19^{2k+1}(24n+5)+1}{24}\right) \equiv 5 \cdot 10^k \cdot e(n) \pmod{19},$$
$$p\left(\frac{19^{2k+2}(24n+23)+1}{24}\right) \equiv 11 \cdot 10^k \cdot f(n) \pmod{19}.$$

*Case where* $m = 23$. Using [Gr-O, Prop. 4], and arguing as above we have for every nonnegative integer $k$

$$F(23, 2k+1; z) \equiv 5^k \eta(24z) E_4(24z) E_6(24z) \pmod{23},$$
$$F(23, 2k+2; z) \equiv 5^{k+1} \eta^{23}(24z) E_4(24z) E_6(24z) \pmod{23}.$$

COROLLARY 12. *Define integers* $g(n)$ *and* $h(n)$ *by*

$$\sum_{n=0}^{\infty} g(n) q^n := E_4(z) E_6(z) \cdot \prod_{n=1}^{\infty} (1-q^n),$$
$$\sum_{n=0}^{\infty} h(n) q^n := E_4(z) E_6(z) \cdot \prod_{n=1}^{\infty} (1-q^n)^{23}.$$

*If $k$ and $n$ are nonnegative integers, then*

$$p\left(\frac{23^{2k+1}(24n+1)+1}{24}\right) \equiv 5^k \cdot g(n) \pmod{23},$$
$$p\left(\frac{23^{2k+2}(24n+23)+1}{24}\right) \equiv 5^{k+1} \cdot h(n) \pmod{23}.$$

PENNSYLVANIA STATE UNIVERSITY, UNIVERSITY PARK, PENNSYLVANIA
*E-mail address*: ono@math.psu.edu